\documentclass[12 pt, a4 paper]{article}
\usepackage{amssymb,amsmath,graphicx,enumerate}
\usepackage{longtable,mathrsfs}
\usepackage{hyperref}
\numberwithin{equation}{section}
\newtheorem{theorem}{Theorem}[section]
\newtheorem{lemma}[theorem]{Lemma}
\newtheorem{definition}{Definition}[section]
\newtheorem{remark}{Remark}[section]

\allowdisplaybreaks

\begin{document}
\begin{center}\begin{large}
On Chebyshev type Inequalities using Generalized k-Fractional Integral Operator
\end{large}\end{center}
\begin{center}
                 $Vaijanath  \, L. Chinchane $

Department of Mathematics,\\
Deogiri Institute of Engineering and Management\\
Studies Aurangabad-431005, INDIA\\
chinchane85@gmail.com
\end{center}
\begin{abstract}
In this paper, using generalized k-fractional integral operator (in terms of the Gauss hypergeometric function), we establish new results on generalized k-fractional integral inequalities by considering the extended Chebyshev functional in case of synchronous function and some other inequalities.
\end{abstract}
\textbf{Keywords :} Chebyshev inequality, generalized k-fractional integral.\\
\textbf{Mathematics Subject Classification :} 26D10, 26A33.
\section{Introduction }
 \paragraph{}In recent years, many authors have worked on fractional integral inequalities by using different fractional integral operator such as Riemann-Liouville, Hadamard, Saigo and Erdelyi-Kober, see \cite{A,BA,BE,C1,C2,C3,C4,C5,D1,KA,P1,YI}. In \cite{KI2}  S. Kilinc and H. Yildirim establish new generalized k-fractional integral inequalities involving Gauss hypergeometric function related to Chebyshev functional. In \cite{C2,D2} authors  gave the following fractional integral inequalities, using the Hadamard and Riemann-Liouville fractional integral for extended Chebyshev functional.
\begin{theorem} Let $f$ and $g$ be two synchronous function on $[0,\infty[$, and $r,p,q:[0,\infty)\rightarrow [0,\infty)$. Then for all $t>0$, $\alpha>0$, we have
 \begin{equation}
 \begin{split}
  &2_{H}D_{1,t}^{-\alpha}r(t) \left[_{H}D_{1,t}^{-\alpha}p(t) _{H}D_{1,t}^{-\alpha}(qfg)(t)+
  _{H}D_{1,t}^{-\alpha}q(t)_{H}D_{1,t}^{-\alpha}(pfg)(t)\right]+\\
  &2 _{H}D_{1,t}^{-\alpha}p(t)_{H}D_{1,t}^{-\alpha}q(t)_{H}D_{1,t}^{-\alpha}(rfg)(t)\geq\\
   &_{H}D_{1,t}^{-\alpha}r(t) \left[_{H}D_{1,t}^{-\alpha}(pf)(t)_{H}D_{1,t}^{-\alpha}(qg)(t)+_{H}D_{1,t}^{-\alpha}(qf)(t)_{H}D_{1,t}^{-\alpha}(pg)(t)\right]+\\
   &_{H}D_{1,t}^{-\alpha}p(t)\left[_{H}D_{1,t}^{-\alpha}(rf)(t)_{H}D_{1,t}^{-\alpha}(qg)(t)+_{H}D_{1,t}^{-\alpha}(qf)(t)_{H}D_{1,t}^{-\alpha}(rg)(t)\right]+\\
   &_{H}D_{1,t}^{-\alpha}q(t)\left[_{H}D_{1,t}^{-\alpha}(rf)(t)_{H}D_{1,t}^{-\alpha}(pg)(t)+_{H}D_{1,t}^{-\alpha}(pf)(t)_{H}D_{1,t}^{-\alpha}(rg)(t)\right]
  \end{split}
\end{equation}
\end{theorem}
\begin{theorem} Let $f$ and $g$ be two synchronous function on $[0,\infty[$, and $r,p,q:[0,\infty)\rightarrow [0,\infty)$. Then for all $t>0$, $\alpha>0$, we have:
 \begin{equation}
 \begin{split}
  &_{H}D_{1,t}^{-\alpha}r(t)\times\\
  & \left[_{H}D_{1,t}^{-\alpha}q(t) _{H}D_{1,t}^{-\beta}(pfg)(t)+2
  _{H}D_{1,t}^{-\alpha}p(t)_{H}D_{1,t}^{-\beta}(qfg)(t)+_{H}D_{1,t}^{-\beta}q(t)_{H}D_{1,t}^{-\alpha}(pfg)(t)\right]\\
  &+\left[_{H}D_{1,t}^{-\alpha}p(t)_{H}D_{1,t}^{-\beta}q(t)+_{H}D_{1,t}^{-\beta}p(t)_{H}D_{1,t}^{-\alpha}q(t)\right]_{H}D_{1,t}^{-\alpha}(rfg)(t)\geq\\
   &_{H}D_{1,t}^{-\alpha}r(t) \left[_{H}D_{1,t}^{-\alpha}(pf)(t)_{H}D_{1,t}^{-\beta}(qg)(t)+_{H}D_{1,t}^{-\beta}(qf)(t)_{H}D_{1,t}^{-\alpha}(pg)(t)\right]+\\
   &_{H}D_{1,t}^{-\alpha}p(t)\left[_{H}D_{1,t}^{-\alpha}(rf)(t)_{H}D_{1,t}^{-\beta}(qg)(t)+_{H}D_{1,t}^{-\beta}(qf)(t)_{H}D_{1,t}^{-\alpha}(rg)(t)\right]+\\
   &_{H}D_{1,t}^{-\alpha}q(t)\left[_{H}D_{1,t}^{-\alpha}(rf)(t)_{H}D_{1,t}^{-\beta}(pg)(t)+_{H}D_{1,t}^{-\beta}(pf)(t)_{H}D_{1,t}^{-\alpha}(rg)(t)\right].
  \end{split}
\end{equation}
\end{theorem}
 \paragraph{}The main objective of this paper is to establish some Chebyshev type inequalities and some other inequalities using generalized k-fractional integral operator. The paper has been organized as follows. In Section 2, we define basic definitions related to generalized k-fractional integral operator. In section 3, we obtain Chebyshev type inequalities using generalized k-fractional. In Section 4 , we prove some inequalities for positive continuous functions.
\section{ Preliminaries}
 \paragraph{}In this section, we present some definitions which will be used later discussion.
\begin{definition}
 Two function $f$ and $g$  are said to synchronous (asynchronous) on $[a,b],$ if
 \begin{equation}
 \left((f(u)-f(v))(g(u)-g(v))\right)\geq (\leq)0,
 \end{equation}
 for all $ u, v \in [0,\infty)$.
 \end{definition}
\begin{definition} \cite{KI2,YI}
The function $f(x)$, for all $x>0$ is said to be in the $L_{p,k}[0,\infty),$ if
\begin{equation}
L_{p,k}[0,\infty)=\left\{f: \|f\|_{L_{p,k}[0,\infty)}=\left(\int_{0}^{\infty}|f(x)|^{p}x^{k}dx\right)^{\frac{1}{p}} < \infty \, \, 1 \leq p < \infty \, k \geq 0\right\},
\end{equation}
\end{definition}
\begin{definition} \cite{KI2,SAO,YI}
Let  $f \in L_{1,k}[0,\infty),$. The generalized Riemann-Liouville fractional integral $I^{\alpha,k}f(x)$ of order $\alpha, k \geq 0$ is defined by
\begin{equation}
I^{\alpha,k}f(x)= \frac{(k+1)^{1-\alpha}}{\Gamma (\alpha)}\int_{0}^{x}(x^{k+1}-t^{k+1})^{\alpha-1}t^{k} f(t)dt.
\end{equation}
\end{definition}
\begin{definition} \cite{KI2,YI}
Let $k\geq0,\alpha>0 \mu >-1$ and $\beta, \eta \in R $. The generalized k-fractional integral $I^{\alpha,\beta,\eta,\mu}_{t,k}$ (in terms of the Gauss hypergeometric function)of order $\alpha$ for real-valued continuous function $f(t)$ is defined by
\begin{equation}
\begin{split}
I^{\alpha,\beta,\eta,\mu}_{t,k}[f(t)]&
= \frac{(k+1)^{\mu+\beta+1}t^{(k+1)(-\alpha-\beta-2\mu)}}{\Gamma (\alpha)}\int_{0}^{t}\tau^{(k+1)\mu}(t^{k+1}-\tau^{k+1})^{\alpha-1}
\times \\
& _{2}F_{1} (\alpha+ \beta+\mu, -\eta; \alpha; 1-(\frac{\tau}{t})^{k+1})\tau^{k} f(\tau)d\tau.
\end{split}
\end{equation}
\end{definition}
where, the function $_{2}F_{1}(-)$ in the right-hand side of (2.4) is the Gaussian hypergeometric function defined by
 \begin{equation}
 _{2}F_{1} (a, b; c; t)=\sum_{n=0}^{\infty}\frac{(a)_{n}(b)_{n}}{(c)_{n}} \frac{t^{n}}{n!},
\end{equation}
and $(a)_{n}$ is the Pochhammer symbol\\
$$(a)_{n}=a(a+1)...(a+n-1)=\frac{\Gamma(a+n)}{\Gamma(a)}, \,\,\,(a)_{0}=1.$$
Consider the function
\begin{equation}
\begin{split}
F(t,\tau)&= \frac{(k+1)^{\mu+\beta+1}t^{(k+1)(-\alpha-\beta-2\mu)}}{\Gamma (\alpha)}\tau^{(k+1)\mu}\\
&(t^{k+1}-\tau^{k+1})^{\alpha-1} \times _{2}F_{1} (\alpha+ \beta+\mu, -\eta; \alpha; 1-(\frac{\tau}{t})^{k+1})\\
&=\sum_{n=0}^{\infty}\frac{(\alpha+\beta+\mu)_{n}(-n)_{n}}{\Gamma(\alpha+n)n!}t^{(k+1)(-\alpha-\beta-2\mu-\eta)}\tau^{(k+1)\mu}(t^{k+1}-\tau^{k+1})^{\alpha-1+n}(k+1)^{\mu+\beta+1}\\
&=\frac{\tau^{(k+1)\mu}(t^{k+1}-\tau^{k+1})^{\alpha-1}(k+1)^{\mu+\beta+1}}{t^{k+1}(\alpha+\beta+2\mu)\Gamma(\alpha)}+\\
&\frac{\tau^{(k+1)\mu}(t^{k+1}-\tau^{k+1})^{\alpha}(k+1)^{\mu+\beta+1}(\alpha+\beta+\mu)(-n)}{t^{k+1}(\alpha+\beta+2\mu+1)\Gamma(\alpha+1)}+\\
&\frac{\tau^{(k+1)\mu}(t^{k+1}-\tau^{k+1})^{\alpha+1}(k+1)^{\mu+\beta+1}(\alpha+\beta+\mu)(\alpha+\beta+\mu+1)(-n)(-n+1)}{t^{k+1}(\alpha+\beta+2\mu+1)\Gamma(\alpha+2)2!}+...
\end{split}
\end{equation}
It is clear that $F(t,\tau)$ is positive  because for all $\tau \in (0, t)$ , $(t>0)$ since each term of the (2.6) is positive.
 \section{Fractional Integral Inequalities for Extended Chebyshev Functional}
 In this section, we establish some Chebyshev type fractional integral inequalities by using the generalized k-fractional integral (in terms of the Gauss hypergeometric function) operator. The following lemma is used for the our main result.
 \begin{lemma} Let $f$ and $g$ be two synchronous function on $[0,\infty[,$ and $x,y:[0,\infty)\rightarrow$ $[0,\infty)$ be two nonnegative functions. Then for all $k \geq 0,$ $t>0$, $\alpha > max\{0,-\beta-\mu\}$, $\beta < 1,$ $\mu >-1,$ $\beta -1< \eta <0,$ we have,
 \begin{equation}
 \begin{split}
  &I^{\alpha,\beta,\eta,\mu}_{t,k}x(t) I^{\alpha,\beta,\eta,\mu}_{t,k}(yfg)(t)+ I^{\alpha,\beta,\eta,\mu}_{t,k}y(t) I^{\alpha,\beta,\eta,\mu}_{t,k}(xfg)(t)\geq \\
  &I^{\alpha,\beta,\eta,\mu}_{t,k}(xf)(t)I^{\alpha,\beta,\eta,\mu}_{t,k}(yg)(t)+I^{\alpha,\beta,\eta,\mu}_{t,k}(yf)(t) I^{\alpha,\beta,\eta,\mu}_{t,k}(xg)(t).
\end{split}
\end{equation}
\end{lemma}
\textbf{Proof}: Since $f$ and $g$  are synchronous on $[0,\infty[$  for all $\tau \geq 0$, $\rho\geq 0$, we have
\begin{equation}
(f(\tau)-f(\rho)) (g(\tau)-g(\rho))\geq 0.
\end{equation}
From (3.2),
\begin{equation}
f(\tau)g(\tau)+f(\rho)g(\rho)\geq f(\tau)g(\rho)+f(\rho)g(\tau).
\end{equation}
\noindent Now, multiplying both side of (3.3) by $ \tau^{k}x(\tau)F(t,\tau)$, $\tau \in (0,t)$, $t>0$. Then the integrating resulting identity with respect to $\tau$  from $0$ to $t$, we obtain by definition (2.4)
\begin{equation}
\begin{split}
&I^{\alpha,\beta,\eta,\mu}_{t,k}(xfg)(t)+f(\rho)g(\rho) I^{\alpha,\beta,\eta,\mu}_{t,k}(x)(t)\\
&I^{\alpha,\beta,\eta,\mu}_{t,k}(yg)(t)I^{\alpha,\beta,\eta,\mu}_{t,k}(xf)(t)+f(\rho)I^{\alpha,\beta,\eta,\mu}_{t,k}(xg)(t).
\end{split}
\end{equation}
\noindent Now, multiplying both side of (3.4) by $ \rho^{k}y(\rho)F(t,\rho)$, $\rho \in (0,t)$, $t>0$, where  $F(t,\rho)$ defined in view of (2.6). Then the integrating resulting identity with respect to $\rho$ from $0$ to $t$, we obtain by definition (2.4)
\begin{equation}
\begin{split}
&I^{\alpha,\beta,\eta,\mu}_{t,k}y(t)I^{\alpha,\beta,\eta,\mu}_{t,k}(xfg)(t)+ I^{\alpha,\beta,\eta,\mu}_{t,k}(yfg)(t)I^{\alpha,\beta,\eta,\mu}_{t,k}(x)(t)\\
&\geq g(\rho)I^{\alpha,\beta,\eta,\mu}_{t,k}(xf)(t)+I^{\alpha,\beta,\eta,\mu}_{t,k}(yf)(t)I^{\alpha,\beta,\eta,\mu}_{t,k}(xg)(t).
\end{split}
\end{equation}
This complete the proof of (3.1)\\
\noindent Now, we gave our  main result here.
\begin{theorem} Let $f$ and $g$ be two synchronous function on $[0,\infty[$, and $r,p,q:[0,\infty)\rightarrow [0,\infty)$. Then for all $k \geq 0,$ $t>0$, $\alpha > max\{0,-\beta-\mu\}$, $\beta < 1,$ $\mu >-1,$ $\beta -1< \eta <0,$ we have,
 \begin{equation}
 \begin{split}
  &2I^{\alpha,\beta,\eta,\mu}_{t,k}r(t) \left[I^{\alpha,\beta,\eta,\mu}_{t,k}p(t)I^{\alpha,\beta,\eta,\mu}_{t,k}(qfg)(t)+
  I^{\alpha,\beta,\eta,\mu}_{t,k}q(t)I^{\alpha,\beta,\eta,\mu}_{t,k}(pfg)(t)\right]+\\
  &2 I^{\alpha,\beta,\eta,\mu}_{t,k}p(t)I^{\alpha,\beta,\eta,\mu}_{t,k}q(t)I^{\alpha,\beta,\eta,\mu}_{t,k}(rfg)(t)\geq\\
   &I^{\alpha,\beta,\eta,\mu}_{t,k}r(t) \left[I^{\alpha,\beta,\eta,\mu}_{t,k}(pf)(t)I^{\alpha,\beta,\eta,\mu}_{t,k}(qg)(t)+I^{\alpha,\beta,\eta,\mu}_{t,k}(qf)(t)I^{\alpha,\beta,\eta,\mu}_{t,k}(pg)(t)\right]+\\
   &I^{\alpha,\beta,\eta,\mu}_{t,k}p(t)\left[I^{\alpha,\beta,\eta,\mu}_{t,k}(rf)(t)I^{\alpha,\beta,\eta,\mu}_{t,k}(qg)(t)+I^{\alpha,\beta,\eta,\mu}_{t,k}(qf)(t)I^{\alpha,\beta,\eta,\mu}_{t,k}(rg)(t)\right]+\\
   &I^{\alpha,\beta,\eta,\mu}_{t,k}q(t)\left[I^{\alpha,\beta,\eta,\mu}_{t,k}(rf)(t)I^{\alpha,\beta,\eta,\mu}_{t,k}(pg)(t)+I^{\alpha,\beta,\eta,\mu}_{t,k}(pf)(t)I^{\alpha,\beta,\eta,\mu}_{t,k}(rg)(t)\right]
  \end{split}
\end{equation}
\end{theorem}
\textbf{Proof}: To prove above theorem, putting $x=p, \ y=q$, and using lemma 3.1, we get
\begin{equation}
\begin{split}
&I^{\alpha,\beta,\eta,\mu}_{t,k}p(t)I^{\alpha,\beta,\eta,\mu}_{t,k}(qfg)(t)+I^{\alpha,\beta,\eta,\mu}_{t,k}q(t)I^{\alpha,\beta,\eta,\mu}_{t,k}(pfg)(t)\geq \\ &I^{\alpha,\beta,\eta,\mu}_{t,k}(pf)(t)I^{\alpha,\beta,\eta,\mu}_{t,k}(qg)(t)+I^{\alpha,\beta,\eta,\mu}_{t,k}(qf)(t)I^{\alpha,\beta,\eta,\mu}_{t,k}(pg)(t).
\end{split}
\end{equation}
\noindent Now, multiplying both side by (3.7) $I^{\alpha,\beta,\eta,\mu}_{t,k}r(t)$, we have
\begin{equation}
\begin{split}
&I^{\alpha,\beta,\eta,\mu}_{t,k}r(t)\left[I^{\alpha,\beta,\eta,\mu}_{t,k}p(t) I^{\alpha,\beta,\eta,\mu}_{t,k}(qfg)(t)+I^{\alpha,\beta,\eta,\mu}_{t,k}q(t)I^{\alpha,\beta,\eta,\mu}_{t,k}(pfg)(t)\right]\geq \\ &I^{\alpha,\beta,\eta,\mu}_{t,k}r(t)\left[I^{\alpha,\beta,\eta,\mu}_{t,k}(pf)(t)I^{\alpha,\beta,\eta,\mu}_{t,k}(qg)(t)+I^{\alpha,\beta,\eta,\mu}_{t,k}(qf)(t)I^{\alpha,\beta,\eta,\mu}_{t,k}(pg)(t)\right],
\end{split}
\end{equation}
\noindent  putting $x=r, y=q$, and using lemma 3.1, we get
 \begin{equation}
 \begin{split}
&I^{\alpha,\beta,\eta,\mu}_{t,k}r(t) I^{\alpha,\beta,\eta,\mu}_{t,k}(qfg)(t)+I^{\alpha,\beta,\eta,\mu}_{t,k}q(t)I^{\alpha,\beta,\eta,\mu}_{t,k}(rfg)(t)\geq \\
&I^{\alpha,\beta,\eta,\mu}_{t,k}(rf)(t)I^{\alpha,\beta,\eta,\mu}_{t,k}(qg)(t)+I^{\alpha,\beta,\eta,\mu}_{t,k}(qf)(t)I^{\alpha,\beta,\eta,\mu}_{t,k}(rg)(t),
\end{split}
\end{equation}
 multiplying both side by (3.9) $I^{\alpha,\beta,\eta,\mu}_{t,k}p(t)$, we have
 \begin{equation}
 \begin{split}
&I^{\alpha,\beta,\eta,\mu}_{t,k}p(t)\left[I^{\alpha,\beta,\eta,\mu}_{t,k}r(t) I^{\alpha,\beta,\eta,\mu}_{t,k}(qfg)(t)+I^{\alpha,\beta,\eta,\mu}_{t,k}q(t)I^{\alpha,\beta,\eta,\mu}_{t,k}(rfg)(t) \right]\geq\\
&I^{\alpha,\beta,\eta,\mu}_{t,k}p(t)\left[I^{\alpha,\beta,\eta,\mu}_{t,k}(rf)(t)I^{\alpha,\beta,\eta,\mu}_{t,k}(qg)(t)+I^{\alpha,\beta,\eta,\mu}_{t,k}(qf)(t)I^{\alpha,\beta,\eta,\mu}_{t,k}(rg)(t)\right].
\end{split}
\end{equation}
With the same arguments as before, we can write
\begin{equation}
\begin{split}
&I^{\alpha,\beta,\eta,\mu}_{t,k}q(t)\left[I^{\alpha,\beta,\eta,\mu}_{t,k}r(t)I^{\alpha,\beta,\eta,\mu}_{t,k}(pfg)(t)+I^{\alpha,\beta,\eta,\mu}_{t,k}p(t)I^{\alpha,\beta,\eta,\mu}_{t,k}(rfg)(t)\right]\geq\\
&I^{\alpha,\beta,\eta,\mu}_{t,k}q(t)\left[I^{\alpha,\beta,\eta,\mu}_{t,k}(rf)(t)I^{\alpha,\beta,\eta,\mu}_{t,k}(pg)(t)+I^{\alpha,\beta,\eta,\mu}_{t,k}(pf)(t)I^{\alpha,\beta,\eta,\mu}_{t,k}(rg)(t)\right].
\end{split}
\end{equation}
Adding the inequalities (3.8), (3.10) and (3.11), we get required inequality (3.6).\\
Here, we give the lemma which is useful to prove our second main result.
\begin{lemma} Let $f$ and $g$ be two synchronous function on $[0,\infty[$. and $x,y:[0,\infty[\rightarrow$ $[0,\infty[$. Then for all $k \geq 0,$ $t>0$, $\alpha > max\{0,-\beta-\mu\}$,$\gamma> max\{0,-\delta-\upsilon\}$ $\beta,\delta < 1,$ $\upsilon,\mu >-1,$ $\beta -1< \eta <0,$ $\delta-1<\zeta <0,$ we have,
 \begin{equation}
 \begin{split}
  &I^{\alpha,\beta,\eta,\mu}_{t,k}x(t) I^{\gamma,\delta,\zeta,\upsilon}_{t,k}(yfg)(t)+ I^{\gamma,\delta,\zeta,\upsilon}_{t,k}y(t) I^{\alpha,\beta,\eta,\mu}_{t,k}(xfg)(t)\geq \\
  &I^{\alpha,\beta,\eta,\mu}_{t,k}(xf)(t) I^{\gamma,\delta,\zeta,\upsilon}_{t,k}(yg)(t)+I^{\gamma,\delta,\zeta,\upsilon}_{t,k}(yf)(t) I^{\alpha,\beta,\eta,\mu}_{t,k}(xg)(t).
\end{split}
\end{equation}
\end{lemma}
\textbf{Proof}:
\noindent Now multiplying both side of (3.4) by
\begin{equation}
\begin{split}
&\frac{(k+1)^{\upsilon+\delta+1}t^{(k+1)(-\delta-\gamma-2\upsilon)}}{\Gamma (\gamma)}\rho^{(k+1)\upsilon}y(\rho)\\
&(t^{k+1}-\rho^{k+1})^{\gamma-1} \times _{2}F_{1} (\gamma+ \delta+\upsilon, -\zeta; \gamma; 1-(\frac{\rho}{t})^{k+1})\rho^{k}
\end{split}
\end{equation}
which remains positive in view of the condition stated in (3.12), $\rho \in (0,t)$, $t>0$, we obtain
\begin{equation}
\begin{split}
&\frac{(k+1)^{\upsilon+\delta+1}t^{(k+1)(-\delta-\gamma-2\upsilon)}}{\Gamma (\gamma)}\rho^{(k+1)\upsilon}y(\rho)\\
&(t^{k+1}-\rho^{k+1})^{\gamma-1} \times _{2}F_{1} (\gamma+ \delta+\upsilon, -\zeta; \gamma; 1-(\frac{\rho}{t})^{k+1})\rho^{k}
I^{\alpha,\beta,\eta,\mu}_{t,k}(xfg)(t)\\
&+\frac{(k+1)^{\upsilon+\delta+1}t^{(k+1)(-\delta-\gamma-2\upsilon)}}{\Gamma (\gamma)}\rho^{(k+1)\upsilon}y(\rho)f(\rho)g(\rho)\\
&(t^{k+1}-\rho^{k+1})^{\gamma-1} \times _{2}F_{1} (\gamma+ \delta+\upsilon, -\zeta; \gamma; 1-(\frac{\rho}{t})^{k+1})\rho^{k}
I^{\alpha,\beta,\eta,\mu}_{t,k}x(t)\geq \\
&\frac{(k+1)^{\upsilon+\delta+1}t^{(k+1)(-\delta-\gamma-2\upsilon)}}{\Gamma (\gamma)}\rho^{(k+1)\upsilon}y(\rho)g(\rho)\\
&(t^{k+1}-\rho^{k+1})^{\gamma-1} \times _{2}F_{1} (\gamma+ \delta+\upsilon, -\zeta; \gamma; 1-(\frac{\rho}{t})^{k+1})\rho^{k}
I^{\alpha,\beta,\eta,\mu}_{t,k}(xf)(t)\\
&+\frac{(k+1)^{\upsilon+\delta+1}t^{(k+1)(-\delta-\gamma-2\upsilon)}}{\Gamma (\gamma)}\rho^{(k+1)\upsilon}y(\rho)f(\rho)\\
&(t^{k+1}-\rho^{k+1})^{\gamma-1} \times _{2}F_{1} (\gamma+ \delta+\upsilon, -\zeta; \gamma; 1-(\frac{\rho}{t})^{k+1})\rho^{k}
I^{\alpha,\beta,\eta,\mu}_{t,k}(xg)(t),
\end{split}
\end{equation}
\noindent then integrating (3.14) over (0,t), we obtain
\begin{equation}
\begin{split}
&I^{\alpha,\beta,\eta,\mu}_{t,k}(xfg)(t)I^{\gamma,\delta,\zeta,\upsilon}_{t,k}y(t)+
I^{\alpha,\beta,\eta,\mu}_{t,k}(x)(t)I^{\gamma,\delta,\zeta,\upsilon}_{t,k}(yfg)(t)\\
 &\geq I^{\alpha,\beta,\eta,\mu}_{t,k}(xf)(t)I^{\gamma,\delta,\zeta,\upsilon}_{t,k}yg(t)
 +I^{\alpha,\beta,\eta,\mu}_{t,k}(xg)(t)I^{\gamma,\delta,\zeta,\upsilon}_{t,k}yf(t),
\end{split}
\end{equation}
\noindent  this ends the proof of inequality (3.12).
\begin{theorem} Let $f$ and $g$ be two synchronous function on $[0,\infty[$, and $r,p,q:[0,\infty)\rightarrow [0,\infty)$. Then for all $t>0$, $\alpha>0$, we have:
 \begin{equation}
 \begin{split}
  &I^{\alpha,\beta,\eta,\mu}_{t,k}r(t)\times\\
  & \left[I^{\alpha,\beta,\eta,\mu}_{t,k}q(t) I^{\gamma,\delta,\zeta,\upsilon}_{t,k}(pfg)(t)+2
  I^{\alpha,\beta,\eta,\mu}_{t,k}p(t)I^{\gamma,\delta,\zeta,\upsilon}_{t,k}(qfg)(t)+I^{\gamma,\delta,\zeta,\upsilon}_{t,k}q(t)I^{\alpha,\beta,\eta,\mu}_{t,k}(pfg)(t)\right]\\
  &+\left[I^{\alpha,\beta,\eta,\mu}_{t,k}p(t)I^{\gamma,\delta,\zeta,\upsilon}_{t,k}q(t)+I^{\gamma,\delta,\zeta,\upsilon}_{t,k}p(t)I^{\alpha,\beta,\eta,\mu}_{t,k}q(t)\right]I^{\alpha,\beta,\eta,\mu}_{t,k}(rfg)(t)\geq\\
   &I^{\alpha,\beta,\eta,\mu}_{t,k}r(t) \left[I^{\alpha,\beta,\eta,\mu}_{t,k}(pf)(t)I^{\gamma,\delta,\zeta,\upsilon}_{t,k}(qg)(t)+I^{\gamma,\delta,\zeta,\upsilon}_{t,k}(qf)(t)I^{\alpha,\beta,\eta,\mu}_{t,k}(pg)(t)\right]+\\
   &I^{\alpha,\beta,\eta,\mu}_{t,k}p(t)\left[I^{\alpha,\beta,\eta,\mu}_{t,k}(rf)(t)I^{\gamma,\delta,\zeta,\upsilon}_{t,k}(qg)(t)+I^{\gamma,\delta,\zeta,\upsilon}_{t,k}(qf)(t)I^{\alpha,\beta,\eta,\mu}_{t,k}(rg)(t)\right]+\\
   &I^{\alpha,\beta,\eta,\mu}_{t,k}q(t)\left[I^{\alpha,\beta,\eta,\mu}_{t,k}(rf)(t)I^{\gamma,\delta,\zeta,\upsilon}_{t,k}(pg)(t)+I^{\gamma,\delta,\zeta,\upsilon}_{t,k}(pf)(t)I^{\alpha,\beta,\eta,\mu}_{t,k}(rg)(t)\right].
  \end{split}
\end{equation}
\end{theorem}
\textbf{Proof}: To prove above theorem, putting $x=p, \ y=q$, and using lemma 3.3 we get
\begin{equation}
\begin{split}
&I^{\alpha,\beta,\eta,\mu}_{t,k}p(t) I^{\gamma,\delta,\zeta,\upsilon}_{t,k}(qfg)(t)+I^{\gamma,\delta,\zeta,\upsilon}_{t,k}q(t)I^{\alpha,\beta,\eta,\mu}_{t,k}(pfg)(t)\geq \\ &I^{\alpha,\beta,\eta,\mu}_{t,k}(pf)(t)I^{\gamma,\delta,\zeta,\upsilon}_{t,k}(qg)(t)+I^{\gamma,\delta,\zeta,\upsilon}_{t,k}(qf)(t)I^{\alpha,\beta,\eta,\mu}_{t,k}(pg)(t).
\end{split}
\end{equation}
\noindent Now, multiplying both side by (3.17) $I^{\alpha,\beta,\eta,\mu}_{t,k}r(t)$, we have
\begin{equation}
\begin{split}
&I^{\alpha,\beta,\eta,\mu}_{t,k}r(t)\left[I^{\alpha,\beta,\eta,\mu}_{t,k}p(t) I^{\gamma,\delta,\zeta,\upsilon}_{t,k}(qfg)(t)+I^{\gamma,\delta,\zeta,\upsilon}_{t,k}q(t)I^{\alpha,\beta,\eta,\mu}_{t,k}(pfg)(t)\right]\geq \\ &I^{\alpha,\beta,\eta,\mu}_{t,k}r(t)\left[I^{\alpha,\beta,\eta,\mu}_{t,k}(pf)(t)I^{\gamma,\delta,\zeta,\upsilon}_{t,k}(qg)(t)+I^{\gamma,\delta,\zeta,\upsilon}_{t,k}(qf)(t)I^{\alpha,\beta,\eta,\mu}_{t,k}(pg)(t)\right],
\end{split}
\end{equation}
\noindent  putting $x=r, \ y=q$, and using lemma 3.3, we get
 \begin{equation}
 \begin{split}
&I^{\alpha,\beta,\eta,\mu}_{t,k}r(t) I^{\gamma,\delta,\zeta,\upsilon}_{t,k}(qfg)(t)+I^{\gamma,\delta,\zeta,\upsilon}_{t,k}q(t)I^{\alpha,\beta,\eta,\mu}_{t,k}(rfg)(t)\geq \\
&I^{\alpha,\beta,\eta,\mu}_{t,k}(rf)(t)I^{\gamma,\delta,\zeta,\upsilon}_{t,k}(qg)(t)+I^{\gamma,\delta,\zeta,\upsilon}_{t,k}qf)(t)I^{\alpha,\beta,\eta,\mu}_{t,k}(rg)(t),
\end{split}
\end{equation}
 multiplying both side by (3.19) $I^{\alpha,\beta,\eta,\mu}_{t,k}p(t)$, we have
 \begin{equation}
 \begin{split}
&I^{\alpha,\beta,\eta,\mu}_{t,k}p(t)\left[I^{\alpha,\beta,\eta,\mu}_{t,k}r(t) I^{\gamma,\delta,\zeta,\upsilon}_{t,k}(qfg)(t)+I^{\gamma,\delta,\zeta,\upsilon}_{t,k}q(t)I^{\alpha,\beta,\eta,\mu}_{t,k}(rfg)(t)\geq \right]\\
&I^{\alpha,\beta,\eta,\mu}_{t,k}p(t)\left[I^{\alpha,\beta,\eta,\mu}_{t,k}(rf)(t)I^{\gamma,\delta,\zeta,\upsilon}_{t,k}(qg)(t)+I^{\gamma,\delta,\zeta,\upsilon}_{t,k}(qf)(t)I^{\alpha,\beta,\eta,\mu}_{t,k}(rg)(t)\right].
\end{split}
\end{equation}
With the same argument as before, we obtain
\begin{equation}
\begin{split}
&I^{\alpha,\beta,\eta,\mu}_{t,k}q(t)\left[I^{\alpha,\beta,\eta,\mu}_{t,k}r(t)I^{\gamma,\delta,\zeta,\upsilon}_{t,k}(pfg)(t)+I^{\gamma,\delta,\zeta,\upsilon}_{t,k}p(t)I^{\alpha,\beta,\eta,\mu}_{t,k}(rfg)(t)\right]\geq \\
&I^{\alpha,\beta,\eta,\mu}_{t,k}q(t)\left[I^{\alpha,\beta,\eta,\mu}_{t,k}(rf)(t)I^{\gamma,\delta,\zeta,\upsilon}_{t,k}(pg)(t)+(pf)(t)I^{\alpha,\beta,\eta,\mu}_{t,k}(rg)(t)\right].
\end{split}
\end{equation}
Adding the inequalities (3.18), (3.20) and (3.21), we follows the inequality (3.16).
\begin{remark}
 If  $ f,g,r,p \  and\  q $ satisfies the following condition,
\begin{enumerate}
  \item The function f and g is asynchronous on $[0,\infty)$.
  \item The function r,p,q are negative on $[0,\infty)$.
  \item Two of the function r,p,q are positive and the third is negative on $[0,\infty)$.
\end{enumerate}
then the inequality 3.6 and 3.16 are reversed.
\end{remark}
\section{Other fractional integral inequalities}
In this section, we proved some fractional integral inequalities for positive and continuous functions which as follows:
\begin{theorem} Suppose that $f$, $g$ and $h$ be three positive and continuous functions on $[0,\infty[$, such that
 \begin{equation}
(f(\tau)-f(\rho))(g(\tau)-g(\tau))(h(\tau)+h(\rho))\geq 0; \  \tau, \rho \in(0,t)\ \  t>0,
\end{equation}
and $x$ be a nonnegative function on $[0,\infty)$. Then for all $k \geq 0,$ $t>0$, $\alpha > max\{0,-\beta-\mu\}$,$\gamma> max\{0,-\delta-\upsilon\}$ $\beta,\delta < 1,$ $\upsilon,\mu >-1,$ $\beta -1< \eta <0,$ $\delta-1<\zeta <0,$ we have,
\begin{equation}
\begin{split}
&I^{\alpha,\beta,\eta,\mu}_{t,k}(x)(t)I^{\gamma,\delta,\zeta,\upsilon}_{t,k}(xfgh)(t)+I^{\alpha,\beta,\eta,\mu}_{t,k}(xh)(t)I^{\gamma,\delta,\zeta,\upsilon}_{t,k}(xfg)(t)\\
&+I^{\alpha,\beta,\eta,\mu}_{t,k}(xfg)(t)I^{\gamma,\delta,\zeta,\upsilon}_{t,k}(xh)(t)+I^{\alpha,\beta,\eta,\mu}_{t,k}(xfgh)(t)I^{\gamma,\delta,\zeta,\upsilon}_{t,k}(x)(t)\\
& \geq I^{\alpha,\beta,\eta,\mu}_{t,k}(xf)(t)I^{\gamma,\delta,\zeta,\upsilon}_{t,k}(xgh)(t)+I^{\alpha,\beta,\eta,\mu}_{t,k}(xg)(t)I^{\gamma,\delta,\zeta,\upsilon}_{t,k}(xfh)(t)\\
&+I^{\alpha,\beta,\eta,\mu}_{t,k}(xgh)(t)I^{\gamma,\delta,\zeta,\upsilon}_{t,k}(xf)(t)+I^{\alpha,\beta,\eta,\mu}_{t,k}(xfh)(t)I^{\gamma,\delta,\zeta,\upsilon}_{t,k}(xg)(t).
\end{split}
\end{equation}
\end{theorem}
\textbf{Proof}: Since $f$, $g$ and $h$ be three  positive and continuous functions on $[0,\infty[$ by (4.1), we can write
\begin{equation}
\begin{split}
&f(\tau)g(\tau)h(\tau)+f(\rho)g(\rho)h(\rho)+f(\tau)g(\tau)h(\rho)+f(\rho)g(\rho)h(\tau)\\
&\geq f(\tau)g(\rho)h(\tau)+f(\tau)g(\rho)h(\rho)+f(\rho)g(\tau)h(\tau)+f(\rho)g(\tau)h(\rho).
\end{split}
\end{equation}
 \noindent Now, multiplying both side of (4.3) by $ \tau^{k}x(\tau)F(t,\tau)$, $\tau \in (0,t)$, $t>0$. Then the integrating resulting identity with respect to $\tau$  from $0$ to $t$, we obtain by definition (2.4)
\begin{equation}
\begin{split}
&I^{\alpha,\beta,\eta,\mu}_{t,k}(xfgh)(t)+f(\rho)g(\rho)h(\rho)I^{\alpha,\beta,\eta,\mu}_{t,k}x(t)+g(\tau)h(\rho)I^{\alpha,\beta,\eta,\mu}_{t,k}(xf)(t)\\
&+f(\rho)g(\rho)I^{\alpha,\beta,\eta,\mu}_{t,k}(xh)(t)\geq g(\rho)I^{\alpha,\beta,\eta,\mu}_{t,k}(xfh)(t)+g(\rho)h(\rho)I^{\alpha,\beta,\eta,\mu}_{t,k}(xf)(t)\\
&+f(\rho)I^{\alpha,\beta,\eta,\mu}_{t,k}(xgh)(t)+f(\rho)h(\rho)I^{\alpha,\beta,\eta,\mu}_{t,k}(xg)(t).
\end{split}
\end{equation}

\noindent Now multiplying both side of (4.4) by
\begin{equation}
\begin{split}
&\frac{(k+1)^{\upsilon+\delta+1}t^{(k+1)(-\delta-\gamma-2\upsilon)}}{\Gamma (\gamma)}\rho^{(k+1)\upsilon}x(\rho)\\
&(t^{k+1}-\rho^{k+1})^{\gamma-1} \times _{2}F_{1} (\gamma+ \delta+\upsilon, -\zeta; \gamma; 1-(\frac{\rho}{t})^{k+1})\rho^{k}
\end{split}
\end{equation}
which remains positive in view of the condition stated in (4.2), $\rho \in (0,t)$, $t>0$ and integrating resulting identity with respective $\rho $ from $0$ to $t$, we obtain
\begin{equation}
\begin{split}
&I^{\alpha,\beta,\eta,\mu}_{t,k}(xfgh)(t)I^{\gamma,\delta,\zeta,\upsilon}_{t,k}x(t)+I^{\gamma,\delta,\zeta,\upsilon}_{t,k}(xfgh)(t)I^{\alpha,\beta,\eta,\mu}_{t,k}x(t)\\
&+I^{\gamma,\delta,\zeta,\upsilon}_{t,k}(xh)(t)I^{\alpha,\beta,\eta,\mu}_{t,k}(xgf)(t)+I^{\gamma,\delta,\zeta,\upsilon}_{t,k}(xfg)(t)I^{\alpha,\beta,\eta,\mu}_{t,k}(xh)(t)\\
&\geq I^{\gamma,\delta,\zeta,\upsilon}_{t,k}xg(t)I^{\alpha,\beta,\eta,\mu}_{t,k}(xfh)(t)+I^{\gamma,\delta,\zeta,\upsilon}_{t,k}(xgh)(t)I^{\alpha,\beta,\eta,\mu}_{t,k}(xf)(t)\\
&+I^{\gamma,\delta,\zeta,\upsilon}_{t,k}(xf)(t)I^{\alpha,\beta,\eta,\mu}_{t,k}(xgh)(t)+I^{\gamma,\delta,\zeta,\upsilon}_{t,k}(xfh)(t)I^{\alpha,\beta,\eta,\mu}_{t,k}(xg)(t).\end{split}
\end{equation}
which implies the proof inequality 4.2.\\
Here, we give another inequality which is as follows.
\begin{theorem}Let $f$, $g$ and $h$ be three positive and continuous functions on $[0,\infty[$, which satisfying the condition (4.1)
 and $x$ and $y$ be two nonnegative functions on $[0,\infty)$. Then for all $k \geq 0,$ $t>0$, $\alpha > max\{0,-\beta-\mu\}$,$\gamma> max\{0,-\delta-\upsilon\}$ $\beta,\delta < 1,$ $\upsilon,\mu >-1,$ $\beta -1< \eta <0,$ $\delta-1<\zeta <0,$ we have,
\begin{equation}
\begin{split}
&I^{\alpha,\beta,\eta,\mu}_{t,k}(x)(t)I^{\gamma,\delta,\zeta,\upsilon}_{t,k}(yfgh)(t)+I^{\alpha,\beta,\eta,\mu}_{t,k}(xh)(t)I^{\gamma,\delta,\zeta,\upsilon}_{t,k}(yfg)(t)\\
&+I^{\alpha,\beta,\eta,\mu}_{t,k}(xfg)(t)I^{\gamma,\delta,\zeta,\upsilon}_{t,k}(yh)(t)+I^{\alpha,\beta,\eta,\mu}_{t,k}(xfgh)(t)I^{\gamma,\delta,\zeta,\upsilon}_{t,k}y(t)\\
& \geq I^{\alpha,\beta,\eta,\mu}_{t,k}(xf)(t)I^{\gamma,\delta,\zeta,\upsilon}_{t,k}(ygh)(t)+I^{\alpha,\beta,\eta,\mu}_{t,k}(xg)(t)I^{\gamma,\delta,\zeta,\upsilon}_{t,k}(yfh)(t)\\
&+I^{\alpha,\beta,\eta,\mu}_{t,k}(xgh)(t)I^{\gamma,\delta,\zeta,\upsilon}_{t,k}(yf)(t)+I^{\alpha,\beta,\eta,\mu}_{t,k}(xfh)(t)I^{\gamma,\delta,\zeta,\upsilon}_{t,k}(yg)(t).
\end{split}
\end{equation}
\end{theorem}
\textbf{Proof}:
 \noindent Multiplying both side of (4.3) by $ \tau^{k}x(\tau)F(t,\tau)$, $\tau \in (0,t)$, $t>0$, where $F(t,\tau)$ defined by (2.6). Then the integrating resulting identity with respect to $\tau$  from $0$ to $t$, we obtain by definition (2.4)
\begin{equation}
\begin{split}
&I^{\alpha,\beta,\eta,\mu}_{t,k}(xfgh)(t)+f(\rho)g(\rho)h(\rho)I^{\alpha,\beta,\eta,\mu}_{t,k}x(t)+g(\tau)h(\rho)I^{\alpha,\beta,\eta,\mu}_{t,k}(xf)(t)\\
&+f(\rho)g(\rho)I^{\alpha,\beta,\eta,\mu}_{t,k}(xh)(t)\geq g(\rho)I^{\alpha,\beta,\eta,\mu}_{t,k}(xfh)(t)+g(\rho)h(\rho)I^{\alpha,\beta,\eta,\mu}_{t,k}(xf)(t)\\
&+f(\rho)I^{\alpha,\beta,\eta,\mu}_{t,k}(xgh)(t)+f(\rho)h(\rho)I^{\alpha,\beta,\eta,\mu}_{t,k}(xg)(t).
\end{split}
\end{equation}

\noindent Now multiplying both side of (4.8) by
\begin{equation}
\begin{split}
&\frac{(k+1)^{\upsilon+\delta+1}t^{(k+1)(-\delta-\gamma-2\upsilon)}}{\Gamma (\gamma)}\rho^{(k+1)\upsilon}y(\rho)\\
&(t^{k+1}-\rho^{k+1})^{\gamma-1} \times _{2}F_{1} (\gamma+ \delta+\upsilon, -\zeta; \gamma; 1-(\frac{\rho}{t})^{k+1})\rho^{k}
\end{split}
\end{equation}
which remains positive in view of the condition stated in (4.7), $\rho \in (0,t)$, $t>0$ and integrating resulting identity with respective $\rho $ from $0$ to $t$, we obtain
\begin{equation}
\begin{split}
&I^{\alpha,\beta,\eta,\mu}_{t,k}(xfgh)(t)I^{\gamma,\delta,\zeta,\upsilon}_{t,k}y(t)+I^{\gamma,\delta,\zeta,\upsilon}_{t,k}(yfgh)(t)I^{\alpha,\beta,\eta,\mu}_{t,k}x(t)\\
&+I^{\gamma,\delta,\zeta,\upsilon}_{t,k}(yh)(t)I^{\alpha,\beta,\eta,\mu}_{t,k}(xgf)(t)+I^{\gamma,\delta,\zeta,\upsilon}_{t,k}(yfg)(t)I^{\alpha,\beta,\eta,\mu}_{t,k}(xh)(t)\\
&\geq I^{\gamma,\delta,\zeta,\upsilon}_{t,k}(yg)(t)I^{\alpha,\beta,\eta,\mu}_{t,k}(xfh)(t)+I^{\gamma,\delta,\zeta,\upsilon}_{t,k}(ygh)(t)I^{\alpha,\beta,\eta,\mu}_{t,k}(xf)(t)\\
&+I^{\gamma,\delta,\zeta,\upsilon}_{t,k}(yf)(t)I^{\alpha,\beta,\eta,\mu}_{t,k}(xgh)(t)+I^{\gamma,\delta,\zeta,\upsilon}_{t,k}(yfh)(t)I^{\alpha,\beta,\eta,\mu}_{t,k}(xg)(t).
\end{split}
\end{equation}
which implies the proof inequality 4.7.

\end{document}